\documentclass[a4paper,11pt,reqno]{article}
\usepackage{amsmath}
\usepackage{amssymb}

\newtheorem{theorem}{Theorem}[section]
\newtheorem{lemma}{Lemma}[section]

\newtheorem{definition}{Definition}[section]
\newtheorem{proposition}{Proposition}[section]

\newtheorem{rem}{Remark}[section]
\def\pf{\noindent{\bf Proof.  }}
\begin{document}
\date{}
\title{\bf  Perturbation analysis for Moore-Penrose inverse of closed operators on Hilbert spaces}
\author{FAPENG DU \thanks{E-mail: jsdfp@163.com}\\
 School of Mathematical and Physical Sciences, Xuzhou Institute of Technology\\
  Xuzhou 221008, Jiangsu Province, P.R. China\\
 \and
 YIFENG XUE \\
 Department of mathematics, East China Normal University\\
  Shanghai 200241, P.R. China
}
\maketitle

\begin{abstract}
In this paper, we investigate the perturbation for the Moore-Penrose inverse
of closed operators on Hilbert spaces. By virtue of a new inner product defined on $H$,
we give the expression of the Moore-Penrose inverse $\bar{T}^\dag$ and the upper bounds of
 $\|\bar{T}^\dag\|$ and $\|\bar{T}^\dag -T^\dag\|$. These results obtained in this paper extend
and improve many related results in this area.


\end{abstract}

\section{Introduction}

An operator $\bar{T}=T+\delta T $ is called the stable perturbation of $T$ if $R(\bar{T}) \cap N(T^+)=\{0\}$.
This notation  is introduced by Chen and the second author in \cite{GX,CWX}.
Later it is generalized to the Banach algebra by the second author in \cite{YFX} and to Hilbert $C^*$--module by
Xu, Wei and Gu in \cite{XWG}. Using this notation the upper bounds for generalized inverse or Moore--Penrose inverse
of bounded linear operators are discussed(See all references).
A classical result about upper bounds is
$$\|\bar{T}^\dag \|\leq \frac{\|T^\dag\|}{1-\|T^\dag\|\|\delta T\|},\quad
\frac{\|\bar{T}^\dag-T^\dag\|}{\|T^\dag\|}\leq\frac{1+\sqrt{5}}{2}\frac{\|T^\dag\|}{1-\|T^\dag\|\|\delta T\|}.$$

In recent years, the perturbation analysis for generalized inverses of closed operators has been
appeared. Some results similar to the perturbation analysis of bounded linear operators are obtained
when $\delta T$ is a T--bounded linear operator(see \cite{QH},\cite{QHL},\cite{WZ}).

But there are some unsolved questions. What is the result of the perturbation for closed operators $T\in C(X,Y)$
when $\delta T$ is a linear operators? What is the expression of the Moore-Penrose inverse
$(T+\delta T)^\dag$ and how to estimate the upper bounds of $\|\bar{T}^\dag\|$ and $\|\bar{T}^\dag -T^\dag\|$
when $X,\,Y$ are Hilbert spaces ? The first question has been solved in \cite{DYX}. Now we discuss the second
question in this paper.

Let $H,K$ be Hilbert spaces, $T\in C(H,K)$ defined on $D(T)$, $\delta T \in L(H,K)$ be a linear operators.
We introduce a new norm $\|\cdot\|_T$ on $D(T)$ such that $(D(T),\|\cdot\|_T)$ be a Hilbert spaces
 and give the expression of $(T+\delta T)^\dag$ and the upper bounds
of $\|\bar{T}^\dag\|$ and $\|\bar{T}^\dag -T^\dag\|$  when $\delta T $ is a bounded linear operators on $(D(T),\|\cdot\|_T)$.

\section{Preliminaries}
Let $X,Y$ be Banach spaces, $L(X,Y),\,C(X,Y)$ and $B(X,Y)$ denote the set of linear operators,
densely-defined closed operators and bounded linear operators from $X$ to $Y$, respectively.
For an operator $T \in L(X,Y)$, $ D(T),\,R(T),\,\ker T$ denoted by the domain, the range and the null spaces of $T$, respectively.

Let $V$ be a closed subspace of $X$. Recall that $V$ is complemented in
$X$ if there is a closed subspace $U$ in $X$ such that $V\cap U=\{0\}$ and $X=V+U$. In this case, we set $X=
V\dotplus U$ and $U=V^c$.

\begin{definition}\rm{\cite{DYX}}
Let $T\in C(X,Y)$. If there is $S\in C(Y,X)$ with $D(S)\supset R(T)$ and $R(S)\subset D(T)$ such that
$$TST=T\; \text{on}\ D(T),\quad STS=S\ \text{on}\; D(S),$$
then $S$ is called a generalized inverse of $T$, which is also denoted by $T^+$.
\end{definition}

Clearly, $P=I-ST$ (resp. $Q=TS$) are idempotent operators on $D(T)$ (resp. $D(S)$)
with $R(P)=\ker T$ (resp. $R(Q)=R(T)$).
\begin{proposition}
Let $T\in C(X,Y)$. Then $T^+\in C(Y,X)$ exists if and only if
$$X=\ker T\oplus \overline{ R(T^+)},\quad Y=\overline{R(T)}\oplus \ker T^+.$$
In addition, $T^+$ is bounded if $R(T)$ closed.
\end{proposition}
\pf
 $(\Rightarrow).$ If $T^+\in C(Y,X)$, then we have
$$D(T)=R(S)+\ker T,\quad D(S)=R(T)+\ker S.$$
So the assertion follows since $D(T)$ (rsep. $D(S)$) are densely in $X$(rsep. $Y$).

$(\Leftarrow).$ See Proposition 2.2 in \cite{DYX}.

\begin{lemma}\label{DL2}\rm{\cite{DYX}}
Let $T\in C(X,Y)$ such that $T^+$ exists. Let $\delta T\colon D(\delta T)\rightarrow D(T^+)$ be a linear operators.
Assume that $I+\delta TT^+\colon D(T^+)\rightarrow D(T^+)$ is bijective.
 Put $\bar{T}=T+\delta T$ and $G=T^+(I+\delta TT^+)^{-1}$.
Then the following statements are equivalent:
\begin{enumerate}
  \item[$(1)$] $R(\bar{T})\cap \ker T^+=\{0\};$
  \item[$(2)$] $\bar{T}G\bar{T}=\bar{T},\; G\bar{T}G=G$ and $R(\bar{T}^+)=R(T^+)$, $\ker \bar{T}^+=\ker T^+$.
  \item[$(3)$] $(I+\delta TT^+)^{-1}\bar{T}$ maps $\ker T$ into $R(T);$
  \item[$(4)$] $(I+\delta T T^+)^{-1}R(\bar{T})=R(T);$
  \item[$(5)$]  $(I+T^+\delta T)^{-1}\ker T=\ker \bar T$.
 \end{enumerate}
\end{lemma}

Let $H$ and $K$ be Hilbert spaces. For $T\in C(H,K)$, let $P_{\overline{R(T)}}$ (resp. $P_{\ker T}$) denote the orthogonal projection
from $K$ (resp. $H$) to $\overline{R(T)}$ (resp. $\ker T$).
\begin{definition}\label{C1}
Let $T\in C(H,K)$. Then there is a unique $S\in C(K,H)$ with $D(S)=R(T)+R(T)^\perp$ and $R(S)=\ker T^\perp\cap D(T)$
such that
\begin{alignat*}{2}
TST&=T\ \text{on}\ D(T),\,\ &\,\ STS&=S\ \text{on}\ D(S),\\
TS&=P_{\overline{R(T)}}\ \text{on}\ D(S),\,\ &\,\ ST&=I-P_{\ker T}\ \text{on}\ D(T).
\end{alignat*}
The operator $S$ is called the Moore--Penrose inverse of $T$, denoted by $T^\dag$. Clearly,
$\ker T^\dag=R(T)^\perp$ and $R(T^\dag)=\ker T^\perp\cap D(T)$.
In addition, if $R(T)$ is closed, then $S$ is bounded.
\end{definition}

\section{Perturbation analysis of M-P inverse on Hilbert spaces}
In this section, we investigate  the expression of M-P inverse $\bar{T}^\dag$ and the upper bound of
$\|\bar{T}^\dag\|$ and $\|\bar{T}^\dag-T^\dag\|$.

$\forall x \in H$, let
$$\|x\|_G=\|x\|+\|Tx\|,$$
then we know $T$ is closed if and only if $(D(T),\|\cdot\|_G)$ is a Banach space(\cite[P191]{TK}).
Clearly $T$ is a bounded linear operators on $(D(T),\|\cdot\|_G)$ since $\|Tx\|\leq \|x\|_G$.

Denote $(\cdot ,\cdot)_H$ be a inner product on $H$.
$\forall x, y \in D(T) $, let $$(x,y)_T=(x,y)_H +(Tx,Ty)_K.$$
It is easy to check that $(x,y)_T$ is a inner product on $D(T)$.

Let $$\|x\|^2_T=(x,x)_T,$$
then $$\|x\|^2_T=(x,x)_T=(x,x)_H +(Tx,Tx)_K=\|x\|^2+\|Tx\|^2,$$
that is, $$\|x\|_T=(\|x\|^2+\|Tx\|^2)^{\frac{1}{2}}.$$
Since
$$ \frac{\sqrt{2}}{2}\|x\|_G \leq \|x\|_T \leq \|x\|_G,$$ we know $\|\cdot\|_G$ equivalence to $\|\cdot\|_T$.
So $T$ is closed if and only if $(D(T),\|\cdot\|_T)$ is a Hilbert space.
 For convenience, we denote $(D(T),\|\cdot\|_T)$ by $D_T$ in the context.

Consider a mapping as following:
\begin{align*}
&\tau :  D(T)\subset H\rightarrow D_T\\
&\tau x=x, \quad \forall x \in D(T)
\end{align*}

Clearly, $\tau $ is defined on $D(T)$ and $R(\tau)=D_T$.

Let $x_n\subset D(T)$ and $x_n\xrightarrow{\|\cdot\|} x,\; \tau x_n\xrightarrow{\|\cdot\|_T} y$, then
$$0\leftarrow \|\tau x_n-y\|^2_T=\|x_n -y\|^2+\|T(x_n -y)\|^2.$$
So $\|x_n -y\|\rightarrow 0$. This indicate $y=\tau x=x \in D(T)$. Hence, $\tau \in C(H,D_T)$.

Clearly,
\begin{align*}
\tau ^\dagger&=\rho \in B(D_T,H);\\
 \rho x&=x,\; x\in D_T.
\end{align*}

\begin{lemma}\rm{\cite{DX}} \label{YLL}
Let $A\in C(L,K),B \in C(H,L)$ with $R(A),R(B), R(AB)$ closed and $R(B)\subseteq D(A)$.
Assume that $AB\in C(H,K)$. Then
\begin{align*}
(AB)^\dag &=P_{\ker (AB)^\perp}(B^\dag (A^\dag ABB^\dag )^\dag A^\dag ) \times \\
&\{A(A^\dag ABB^\dag )(A^\dag ABB^\dag )^\dag A^\dag  +(A^\dag )^*(A^\dag ABB^\dag )(A^\dag ABB^\dag )^\dag A^*-I\}^{-1}.
\end{align*}
\end{lemma}

\begin{lemma}\label{TYL}
Let $T \in C(H,K)$, then $T^+\in B(K,H)$ if and only if $T^+ \in B(K,D_T)$,
and in  this case $$\|T^+\|^2 \leq \|T^+\|^2_T \leq \|T^+\|^2+\|TT^+\|^2.$$
\end{lemma}
\pf

If $T^+ \in B(K,H)$, then  $TT^+ \in B(K)$. $\forall x \in K$,
$$\|T^+x\|^2_T=\|T^+x\|^2+\|T(T^+x)\|^2\leq (\|T^+\|^2+\|TT^+\|^2)\|x\|^2.$$
Hence, $T^+ \in B(K,D_T)$ and $\|T^+\|^2_T \leq \|T^+\|^2+\|TT^+\|^2$.

Conversely, if $T^+ \in B(K,D_T)$, then $\forall x \in K$,
$$\|T^+x\|^2=\|T^+x\|^2_T-\|TT^+x\|^2 \leq \|T^+x\|^2_T.$$
Hence,  $T^+ \in B(K,H)$ and $\|T^+\|\leq \|T^+\|_T$.

From the above, we have
$$\|T^+\|^2 \leq \|T^+\|^2_T \leq \|T^+\|^2+\|TT^+\|^2.$$

\begin{lemma}\label{YL2}
Let $T\in C(H,K)$ with $R(T)$ closed. If $T$ has generalized inverse $T^+$,
then $T^\dag \in B(K,H)$ and
$$T^\dag=-P_{\ker T^\perp}(I+P(I-P-P^*)^{-1})T^+(I-Q-Q^*)^{-1}.$$
\end{lemma}
\pf Since $R(T)$ closed, we have $T^+\in B(K,H)$. So $T^+ \in B(K,D_T)$ by Lemma \ref{TYL}.
Thus, $Q=TT^+ \in B(K),\; P=I-T^+T \in B(D_T)$ are idempotent operators. Now we consider the Moore-Penrose inverse $T^\dag$ of
$T$ on $D_T$. From \cite{WXQ}, we have $T^\dag \in B(K,D_T)$ and
$$T^\dag=-(I+P(I-P-P^*)^{-1})T^+(I-Q-Q^*)^{-1}.$$
Since $T^\dag \in B(K,D_T)$, we have $T^\dag \in B(K,H)$ by Lemma \ref{TYL}.
Noting that $T\in C(H,K)$ is a compound operator by $T\in B(D_T,K)$ and $\tau \in C(H,D_T)$.
Therefore, by Lemma \ref{YLL}, we have
$$T^\dag=-P_{\ker T^\perp}(I+P(I-P-P^*)^{-1})T^+(I-Q-Q^*)^{-1}.$$

\begin{theorem}\label{dl2}
Let $T\in C(H,K)$ with $T^\dag \in B(K,H)$, $\delta T \in B(D_T,K)$
 such that $\bar{T}=T+\delta T$ closed, $D(T)\subseteq D(\delta T)$. If
 $I+\delta T T^\dag$ is invertible and $R(\bar{T})\cap N(T^\dag)=\{0\}$,
then $\bar{T}^\dag \in B(K,H)$ and
$$\bar{T}^\dag=-P_{\ker \bar{T}^\perp}(I+\bar{P}(I-\bar{P}-\bar{P}^*)^{-1})G(I-\bar{T}G-(\bar{T}G)^*)^{-1},$$
where $G=T^\dag(I+\delta TT^\dag)^{-1},\;\bar{P}=I-G\bar{T}$.
\end{theorem}
\pf  $\forall x \in D(T)$, there is an $M$ such that
$\|\delta Tx\|\leq M\|x\|_T$ since $\delta T\in B(D_T, K)$. Thus, $\forall y \in K$,
$$\|\delta T T^\dag y\|^2\leq \|\delta T\|^2_T \|T^\dag y\|^2_T\leq \|\delta T\|^2_T (\|T^\dag\|^2+1)\|y\|^2.$$
Hence $G =T^\dag(I+\delta TT^\dag)^{-1} \in B(K,H)$ be the generalized inverse of $\bar{T}$ by Lemma \ref{DL2}.
By Lemma \ref{YL2}, $\bar{T}^\dag \in B(K,H)$ and
$$\bar{T}^\dag=-P_{\ker \bar{T}^\perp}(I+\bar{P}(I-\bar{P}-\bar{P}^*)^{-1})G(I-\bar{T}G-(\bar{T}G)^*)^{-1}.$$

\begin{rem}
If $\delta T$ is T--bounded, i.e.,
there are constants $a,\,b>0$ such that
$$\|\delta Tx\|\le a\|x\|+b\|Tx\|,\quad \forall\,x\in D(T),$$
then $\delta T \in B(D_T,K)$. Indeed,
$$\|\delta Tx\|^2\leq (a\|x\|+b\|Tx\|)^2\leq 2(\max(a,b))^2(\|x\|^2+\|Tx\|^2)=2(\max(a,b))^2\|x\|^2_T.$$
\end{rem}

Let $M,N$ are two closed subspaces of $H$. Set
$$\delta(M,N)=\sup\{dist(\mu,N)|\|\mu\|=1,\mu \in M \}.$$
We call
$\hat{\delta}(M,N)=\max\{\delta(M,N),\delta(N,M)\}$
the gap between subspaces $M$ and $N$.
\begin{proposition}\rm{\cite{TK}}
\begin{description}
  \item[(1)] $\delta(M,N)=0$ if and only if $M\subset N$
  \item[(2)] $\hat{\delta}(M,N)=0$ if and only if $M=N$
  \item[(3)] $\hat{\delta}(M,N)=\hat{\delta}(N,M)$
  \item[(4)] $0\leq \delta(M,N)\leq 1$, $0\leq \hat{\delta}(M,N)\leq 1$
  \item[(5)]$\hat{\delta}(M,N)=\|P-Q\|$, where $P,Q$ are orthogonal projection on $M,N$, respectively.
\end{description}
\end{proposition}
For convenience, we set
$\|\delta T\|_T=\underset{\|x\|_T\leq 1}{\sup}\dfrac{\|\delta Tx\|}{\|x\|_T}$ if $\delta T \in B(D_T,K)$.

\begin{lemma}
Under the assumptions of Theorem \ref{dl2}, we have
\begin{enumerate}
  \item $\delta(R(T),R(\bar{T}))\leq \|\delta T\|_T(\|T^\dag\|^2+1)^\frac{1}{2}.$
  \item $\delta(\ker T,\ker (\bar{T}))\leq \|\bar{T}^\dag\|\|\delta T\|_T.$
\end{enumerate}
\end{lemma}
\pf  Noting that $\|\delta T x\|\leq \|\delta T\|_T\|x\|_T$ and $\bar{T}^\dag \in B(K,H)$.

$(1).$  Let $u\in R(T)$ with $\|u\|=1$, then there is a $x\in D(T)$ such that $u=Tx$.
\begin{align*}
dist(u,R(\bar{T}))&\leq \|u-\bar{T}(T^\dag Tx)\|=\|\delta TT^\dag Tx\| \\
&\leq \|\delta T\|_T\|T^\dag Tx\|_T \\
&\leq \|\delta T\|_T(\|T^\dag \|^2+1)^\frac{1}{2}\|u\|.
\end{align*}
Hence, $\delta(R(T),R(\bar{T}))\leq \|\delta T\|_T(\|T^\dag \|^2+1)^\frac{1}{2}$.

$(2).$
Let $x\in \ker T$ with $\|x\|=1$, then $Tx=0$
\begin{align*}
dist(x,\ker (\bar{T}))&\leq \|x-(I-\bar{T}^\dag \bar{T})x\|=\|\bar{T}^\dag \delta Tx\| \\
&\leq \|\bar{T}^\dag\|\|\delta Tx\| \\
&\leq\|\bar{T}^\dag\|_T\|\delta T\|_T\|x\|.
\end{align*}
Hence, $\delta(\ker T,\ker (\bar{T}))\leq \|\bar{T}^\dag\|\|\delta T\|_T$.

\begin{lemma}\label{TY3}
Let M and N be closed subspaces of H. Suppose that $M\cap N^\perp=\{0\}$.Then
$$\delta (M,N)=\|P_M-P_N\|.$$
\end{lemma}
\pf If $\delta (M,N)=1$, then $1=\hat{\delta}(M,N)=\|P_M-Q_N\|$, Thus $\delta (M,N)=\|P_M-P_N\|$.\\
Assume that $\delta (M,N)=\delta <1$, then $\forall x \in M$,
$$ \|(I-P_N)P_Mx\|=dist(P_Mx,N)\leq \|P_Mx\|\delta (M,N)\leq \delta \|x\|.$$
So by Lemma 3 of \cite{XCC}, we know $\delta (M,N)=\|P_M-P_N\|$.

\begin{lemma}
Under the assumptions of Theorem \ref{dl2}, we have
$$ \|TT^\dag -\bar{T}\bar{T}^\dag\|= \delta(R(T),R(\bar{T})).$$
\end{lemma}
\pf Since $T^\dag \in B(K,H)$, we have $\ker (T^\dag)=R(T)^\perp$. So $\ker (T^\dag)=\ker (\bar{T}^\dag)$
implies $R(T)\cap \ker (\bar{T}^\dag)=\{0\}$.
By Lemma  \ref{TY3}, we know $$ \|TT^\dag -\bar{T}\bar{T}^\dag\|= \delta(R(T),R(\bar{T})).$$

\begin{theorem}\label{Dl3}
Under the assumptions of Theorem \ref{dl2}, we have
\begin{description}
\item[(1)] $\|\bar{T}^\dag\|\leq \|(1+\delta T T^\dag)^{-1}\|(\|T^\dag\|^2+1)^\frac{1}{2}$.
\item[(2)] $\|\bar{T}^\dag-T^\dag\|\leq \dfrac{1+\sqrt{5}}{2}\|\bar{T}^\dag\|\|\delta T\|_T(\|T^\dag\|^2+1)^\frac{1}{2}$.
\end{description}

\end{theorem}
\pf
$(1).$ Since $\delta T \in B(D_T,K)$, we have
$$\|\delta T T^\dag x\|\leq \|\delta T\|_T \|\|T^\dag x\|_T \leq \|\delta T\|_T (\|T^\dag \|^2+1)^\frac{1}{2}\|x\|.$$
Hence, $\delta TT^\dag \in B(K)$. By Lemma \ref{DL2}  ,we have
$$G=T^+(I+\delta T T^+)^{-1} \in B(K,H),\; P=I-G\bar{T}\in B(D_T),\;Q=\bar{T}G \in B(K)$$
and $\|(I+P(I-P-P^*)^{-1})x\|_T\leq \|x\|_T,\; x\in D(T).$ So,
\begin{align*}
\|\bar{T}^\dag x\|^2&\leq\|\bar{T}^\dag x\|^2_T=\|-(I+P(I-P-P^*)^{-1})G(I-Q-Q^*)^{-1}x\|^2_T \\
&\leq\|G(I-Q-Q^*)^{-1}x\|^2_T\\
&\leq \|T^\dag(1+\delta T T^\dag)^{-1}(I-Q-Q^*)^{-1}x\|^2+\|TT^\dag(1+\delta T T^\dag)^{-1}(I-Q-Q^*)^{-1}x\|^2\\
&\leq (\|T^\dag\|^2+1)\|(1+\delta T T^\dag)^{-1}\|^2\|x\|^2.
\end{align*}
$(2).$
Since  $D(T)$ dense in $H$, we can extend $I-T^\dag T$ to the whole spaces $H$ such that $P_{\ker T}|_{D(T)}=I-T^\dag T$
and $P_{\ker T}$ is an orthogonal projection form $H$ onto $\ker T$. Similarly, we extend $I-\bar{T}^\dag \bar{T}$
to the whole spaces $H$ such that $P_{\ker (\bar{T})}|_{D(\bar{T})}=I-\bar{T}^\dag \bar{T}$
and $P_{\ker (\bar{T})}$ is an orthogonal projection form $H$ onto $\ker (\bar{T})$.

Clearly, $\forall x \in D(T)$, $(T^\dag T-\bar{T}^\dag \bar{T})x=(P_{\ker (\bar{T})}-P_{\ker T})x.$

Noting that $\ker T\cap \ker(\bar{T})^\perp=\{0\}$, we have
$\|P_{\ker T}-P_{\ker (\bar{T})}\|= \delta(\ker T,\ker (\bar{T}))$ by Lemma \ref{TY3}.

$\forall y \in K$ and $\|y\|=1$,
\begin{align*}
\bar{T}^\dag y-T^\dag y
&=-\bar{T}^\dag \delta TT^\dag TT^\dag y+\bar{T}^\dag (\bar{T}\bar{T}^\dag -TT^\dag )(I-TT^\dag )y\\
&+(P_{\ker (\bar{T})}-P_{\ker T})T^\dag y.
\end{align*}
Using the proof of Proposition 7 in \cite{XCC}, we have
$$\|\bar{T}^\dag y-T^\dag y \|^2_T \leq \frac{3+\sqrt{5}}{2}\|\bar{T}^\dag \|^2\|\delta T\|^2_T(\|T^\dag \|^2+1).$$
Since $ \|\bar{T}^\dag y-T^\dag y\|\leq \|\bar{T}^\dag y-T^\dag y\|_T $,
we have
$$\|\bar{T}^\dag-T^\dag\|\leq \frac{1+\sqrt{5}}{2}\|\bar{T}^\dag \|\|\delta T\|_T(\|T^\dag \|^2+1)^\frac{1}{2}.$$

\section{Perturbation analysis for $Tx=b$ in Hilbert spaces}

In this section, we consider the perturbation of the least square solution
of the following two equations
\begin{description}
  \item[(1)] $Tx=b, $
  \item[(2)] $\bar{T}\bar{x}=\bar{b},$ $(\bar{b}=b+\delta b)$
\end{description}

 As we know the solutions of $$\|Tx-b\|=\min_{z\in D(T)} \|Tz-b\|$$ are
$x=T^\dag b+(I-T^\dag T)z,\forall z\in D(T)$, denoted by $S(T,b)$,
i.e. $$S(T,b)=\{x :x=T^\dag b+(I-T^\dag T)z,\forall z\in D(T)\}$$
Similarly,
$$S(\bar{T},\bar{b})=\{\bar{x}:\bar{x}=\bar{T}^\dag \bar{b}+(I-\bar{T}^\dag \bar{T})z,\forall z\in D(\bar{T})\}$$
\begin{theorem}
 Under the assumptions of Theorem \ref{dl2}, we have\\
$(1)$ For any solution
$x=T^\dag b+(I-T^\dag T)z$ in $S(T,b)$,
there exist  $\bar{x}\in S(\bar{T},\bar{b})$ such that
$$\|\bar{x}- x\|\leq \|\bar{T}^\dag\|\|(b-Tx)+(\delta b-\delta Tx)\|.$$
$(2)$ For any solution $\bar{x}=\bar{T}^\dag \bar{b}+(I-\bar{T}^\dag \bar{T})z$ in $S(\bar{T},\bar{b})$,
there exist  $x\in S(T,b)$ such that
$$\|\bar{x}- x\|\leq \|T^\dag\|\|(\bar{b}-\bar{T}\bar{x})-(\delta b-\delta T\bar{x})\|.$$
\end{theorem}
\pf $(1)$
Taking
$$\bar{x}=\bar{T}^\dag \bar{b}+(I-\bar{T}^\dag \bar{T})(T^\dag b+(I-T^\dag T)z).$$
Then
\begin{align*}
\|\bar{x}- x\|&=\|\bar{T}^\dag \bar{b}+(I-\bar{T}^\dag \bar{T})(T^\dag b+(I-T^\dag T)z)- (T^\dag b+(I-T^\dag T)z)\|\\
&=\|\bar{T}^\dag \bar{b}-(\bar{T}^\dag \bar{T})(T^\dag b+(I-T^\dag T)z)\|\\
&=\|\bar{T}^\dag \delta b+(I-\bar{T}^\dag \bar{T})T^\dag b -\bar{T}^\dag \bar{T}(I-T^\dag T)z +(\bar{T}^\dag-T^\dag)b\| \\
&=\|\bar{T}^\dag \delta b+(I-\bar{T}^\dag \bar{T})T^\dag b -\bar{T}^\dag \bar{T}(I-T^\dag T)z \\
&+(-\bar{T}^\dag \delta TT^\dag b+\bar{T}^\dag (I-TT^\dag)b-(I-\bar{T}^\dag \bar{T})T^\dag b)\|\\
&=\|\bar{T}^\dag(\bar{b}-\bar{T}x)\|\\
&\leq \|\bar{T}^\dag\|\|(b-Tx)+(\delta b-\delta Tx)\|.
\end{align*}
$(2)$
Taking
$$x=T^\dag b+(I-T^\dag T)(\bar{T}^\dag \bar{b}+(I-\bar{T}^\dag \bar{T})z). $$
By similarly computation to $(1)$, we have
$$\|\bar{x}- x\|\leq \|T^\dag\|\|(\bar{b}-\bar{T}\bar{x})-(\delta b-\delta T\bar{x})\|.$$

\begin{theorem}
 Under the assumptions of Theorem \ref{dl2}, we have
\begin{align*}
\|\bar{T}^\dag \bar{b}-T^\dag b\|&\leq\|(1+\delta T T^\dag)^{-1}\|(\|T^\dag\|^2+1)^\frac{1}{2}\|\delta b\|\\
&+ \frac{1+\sqrt{5}}{2}\|(1+\delta T T^\dag)^{-1}\|\|\delta T\|_T(1+\|T^\dag \|^2)^{\frac{1}{2}}\|b\|.
\end{align*}
\end{theorem}
\pf By Theorem \ref{Dl3}, we have
\begin{align*}
\|\bar{T}^\dag \bar{b}-T^\dag b\|&=\|\bar{T}^\dag \delta b +(\bar{T}^\dag -T^\dag )b\|\\
&\leq \|\bar{T}^\dag \delta b\|+\|\bar{T}^\dag -T^\dag \|\|b\|\\
&\leq \|(1+\delta T T^\dag)^{-1}\|(\|T^\dag\|^2+1)^\frac{1}{2}\|\delta b\|\\
&+ \frac{1+\sqrt{5}}{2}\|(1+\delta T T^\dag)^{-1}\|\|\delta T\|_T(1+\|T^\dag \|^2)^{\frac{1}{2}}\|b\|.
\end{align*}

\section{Conclusion}
In this paper, we extend the perturbation analysis of Moore-Penrose inverse for bounded linear operators to closed operators.
By virtue of a new inner product, we give the expression of the Moore-Penrose inverse $\bar{T}^\dag$ and the upper bounds of
 $\|\bar{T}^\dag\|$ and $\|\bar{T}^\dag -T^\dag\|$. As an application, we study the perturbation of the least square solution.
 These results enrich and improve the perturbation theory of Moore-Penrose
 inverse described in \cite{Xue}.

\end{document}